\documentclass[11pt]{article}
\usepackage{amsmath}
\usepackage{amssymb}
\usepackage{latexsym}
\usepackage{epsf}

\usepackage[cp866]{inputenc}
\usepackage[english]{babel}

\usepackage{amsthm, amstext}
\usepackage{array, amsfonts, mathrsfs}

\theoremstyle{plain}
\newtheorem{T}{Theorem}

\newcommand{\re}[1]{(\ref{#1})}

\title{Asymptotic behavior of singular values of the acoustic observation problem}
\author{M. N. Demchenko\footnote{St.Petersburg Department of V.A.~Steklov Institute of Mathematics of Russian Academy of Sciences. 
		demchenko@pdmi.ras.ru.}}

\date{}

\begin{document}

\maketitle

\begin{abstract}
We consider the problem of recovering of initial data in the IBVP for the wave-type equation in the half-space by the solution restricted to the boundary. The singular value decomposition of this problem is concerned: the asymptotics of singular values is obtained. 
\smallskip

\end{abstract}

\medskip

\section{Introduction}

Fix an integer $d \geqslant 2$ and put ${\mathbb R}^d_+ = \{x=(x_1, \cdots, x_d)\in {\mathbb R}^d\,|\, x_d > 0\}$.
Suppose $q \in C_0^\infty({\mathbb R}^d_+)$ is a real-valued function. 
Consider the following initial boundary value problem for scalar function $u(x,t)$, $(x,t)\in {\mathbb R}^d_+\times {\mathbb R}$:
\begin{align}
	&\partial^2_t u - \Delta u + q u = 0, \quad
	\partial_{x_d} u|_{x_d = 0} = 0, \notag\\
	&u|_{t=0} = 0, \quad \partial_t u|_{t=0} = v,
	\label{BVP}
\end{align}
where $v \in C_0^\infty({\mathbb R}^d_+)$.
Introduce the operator 
\[
	{\cal O}: v \mapsto u|_{\Sigma_0},
\]
where $\Sigma_0 =  \partial{\mathbb R}^d_+\times {\mathbb R}$ (the time-space boundary) 
and $u$ is the solution of \re{BVP} for the given initial data $v$.
The operator ${\cal O}$ is well defined as the solution $u$ is regular.

Next we introduce a restriction of the operator ${\cal O}$.
Let $\Omega \subset {\mathbb R}^d_+$ be a bounded open set that satisfies $\overline\Omega \subset {\mathbb R}^d_+$.
Also let $\Sigma$ be a bounded relatively open subset of $\Sigma_0$, such that $\overline\Sigma\subset\partial{\mathbb R}^d_+\times (0,\infty)$.
Consider the operator ${\cal O}_\Omega^\Sigma: C_0^\infty(\Omega) \to C^\infty(\Sigma)$ acting as follows
\[
	{\cal O}_\Omega^\Sigma \, v = ({\cal O} v)|_{\Sigma}.
\]
The operator ${\cal O}_\Omega^\Sigma$ can be continued as compact operator acting from $L_2(\Omega)$ to $L_2(\Sigma)$ 
(see sec.~\ref{secProof}).

We consider the problem of recovering $v$ by ${\cal O}_\Omega^\Sigma v$ (the coefficient $q$ is given), which we call the {\em observation problem} after~\cite{BLR}.
Invertibility of ${\cal O}_\Omega^\Sigma$ depends on the geometry of $\Omega$ and $\Sigma$.
An example of the case when ${\cal O}_\Omega^\Sigma$ is invertible (more precisely ${\rm Ker}\, {\cal O}_\Omega^\Sigma = \{0\}$)
give $\Sigma = \Gamma\times(0,T)$ and $\Omega$ satisfying $\Omega \subset {\mathbb R}^d_+ \cap (\cup_{\gamma\in\Gamma} B_\gamma(T))$. Here $T>0$, $\Gamma$ is an open bounded subset of $\partial{\mathbb R}^d_+$, $B_\gamma(T)$ is an open ball of radius $T$ centered in $\gamma$ \cite{Bel}.
Note that similar problems were also considered in~\cite{BuhKlib, Klib}, where the method based on Carleman estimates was developed.
In this paper we are interested in the procedure of recovering $v$ based on the
singular value decomposition (SVD) for ${\cal O}_\Omega^\Sigma$.
Recall that for a compact linear operator $A$ acting in Hilbert spaces $H_0 \to H_1$ there exist orthonormal basis $\{v_n\} \subset H_0$ and $\{f_n\} \subset H_1$ such that
\[
	A v = \sum_{n\geqslant 1} s_n(A)\,  (v, v_n)_{H_0}\, f_n . 
\]
Here $s_n(A)$ are singular values of $A$ defined as $s_n(A) = (\lambda_n(A^*A))^{1/2}$, where $\lambda_n(A^*A)$ are eigenvalues of the compact operator $A^* A$ being numbered in non-increasing order with the multiplicity taken into account.
If ${\rm Ker}\, A = \{0\}$ and $f$ belongs to the range of $A$ then
\begin{equation}
	A^{-1} f = \sum_{n\geqslant 1} s_n(A)^{-1}\,  (f, f_n)_{H_1}\, v_n.
	\label{series}
\end{equation}
It can be seen that the stability of inversion 
(in particular) depends on the behavior of $s_n(A)$. 
In this paper we prove that 
$s_n({\cal O}_\Omega^\Sigma)$ decrease as a negative power of $n$ (relation~\re{asymp_sn}). Note that we do not need the invertibility of ${\cal O}_\Omega^\Sigma$ and thus no geometric conditions on $\Omega$ and $\Sigma$ are imposed.

To study the asymptotics of $s_n({\cal O}_\Omega^\Sigma)$ we approximate the operator ${\cal O}_\Omega^\Sigma$ 
by Fourier integral operators (FIO), which yields approximation of $({\cal O}_\Omega^\Sigma)^* {\cal O}_\Omega^\Sigma$ by pseudodifferential operators ($\Psi$DO) of order $-2$. Then we apply the result on 
asymptotics of singular values for non-elliptic $\Psi$DOs of negative order \cite{BS_anisotropic}.
The theory of FIOs is a natural tool for analysis of operator ${\cal O}_\Omega^\Sigma$.
Furthermore, FIOs were applied
to the problems of integral geometry 
\cite{FrikQuinto, GuS, Nguyen, Palamodov_Remarks, Uhlmann}
(the list of references is not complete),
which are similar to the problem of inverting ${\cal O}_\Omega^\Sigma$.

The principal symbols of $\Psi$DOs that approximate $({\cal O}_\Omega^\Sigma)^* {\cal O}_\Omega^\Sigma$ vanish on the subset of $T^* \Omega \setminus 0$ of nonzero measure.
If $({\cal O}_\Omega^\Sigma)^* {\cal O}_\Omega^\Sigma$ was precisely an elliptic $\Psi$DO (i.e. operator with nonvanishing principal symbol) then $(({\cal O}_\Omega^\Sigma)^* {\cal O}_\Omega^\Sigma)^{-1}$ would be a $\Psi$DO of order $2$; if additionally ${\rm Ker}\, {\cal O}_\Omega^\Sigma = \{0\}$ then the operator $({\cal O}_\Omega^\Sigma)^{-1}$ would be bounded in certain Sobolev spaces since
\[
	({\cal O}_\Omega^\Sigma)^{-1} = (({\cal O}_\Omega^\Sigma)^* {\cal O}_\Omega^\Sigma)^{-1} ({\cal O}_\Omega^\Sigma)^*.
\]
In our case the operators that approximate $({\cal O}_\Omega^\Sigma)^* {\cal O}_\Omega^\Sigma$ do not have even parametrices, so 
the observation problem is ill-posed.
The principal symbols of approximating $\Psi$DOs vanish outside of the
conic set $AZ \subset T^*\Omega\setminus 0$ which is defined in sec.~\ref{mainres}.
The same situation arises in problems of integral geometry with {\em limited data}~-- see 
\cite{Begmatov, FrikQuinto, Nguyen, Palamodov, Uhlmann}.
The conic set $AZ$ is called an {\em audible zone} (or {\em visible zone}) there.
The results of these papers show that it is possible to reconstruct some part of singularities of unknown function, namely~-- the intersection of the wavefront set of the function with the audible zone. 
This gives a hope that series~\re{series} applied to ${\cal O}_\Omega^\Sigma$ allows to reconstruct $v$ microlocally in the audible zone, the reconstruction being stable in some reasonable sense.
However, in the present paper we only obtain the asymptotics of $s_n({\cal O}_\Omega^\Sigma)$.

Note that for some problems of integral geometry singular values and singular basis (i.e. $v_n$ and $f_n$ in~\re{series}) were found in a nearly explicit form.
It concerns some particular (though very important) cases such as limited angle tomography \cite{Louis, Natterer} and exterior Radon transform \cite{QuintoSVD}.
It seems that there are no such results for ${\cal O}_\Omega^\Sigma$ even for $q=0$ and some certain $\Omega$ and $\Sigma$.

\textbf{Acknowledgements.}
The research was supported by RFBR grants 14-01-31388 mol-a, 15-31-20600 mol-a-ved, 14-01-00535 and SPbGU 6.38.670.2013.

The author thanks prof. M.I. Belishev for helpful discussions on this problem.

\section{The audible zone and the asymptotics of $s_n({\cal O}_\Omega^\Sigma)$} 
\label{mainres}
To define the audible zone $AZ$ and formulate the main result we introduce some notation.
A point in the plane $\Sigma_0$ will be identified by $(x', t)$, $x' \in {\mathbb R}^{d-1}$, $t\in {\mathbb R}$ (we use identification of $\partial{\mathbb R}^d_+$ and ${\mathbb R}^{d-1}$).
For $x \in {\mathbb R}^d$ we put $x' = (x_1, \dots, x_{d-1})$.
We also use natural identification of cotangent bundle $T^*\Omega$ and $\Omega\times{\mathbb R}^d$. 
Introduce the following mapping from $T^*\Omega$ to $\Sigma_0$:
\[
	\gamma(y, \eta) = \left(y' - y_d\, \frac{\eta'}{\eta_d},\, y_d \, \frac{|\eta|}{|\eta_d|}\right),
	\quad y \in \Omega,\, \eta\in T^*_y\Omega
\]
(the function $\gamma$ is defined if $\eta_d \ne 0$, i.e. almost everywhere in $T^*\Omega$).
Note that the first component of $\gamma(y,\eta)$ is an intersection point of the line $\{y + s\eta,\, s\in{\mathbb R}\}$ and $\partial{\mathbb R}^d_+$, while the second component is the distance between this intersection point and $y$.
The audible zone is defined as follows
\[
	AZ = \{(y,\eta)\in T^*\Omega\setminus 0\,|\,\,\eta_d\ne 0,\, \gamma(y,\eta)\in\Sigma\}.
\]
We say that $K\subset T^*\Omega\setminus 0$ is a {\em conic set} if it is invariant with respect to the mapping $(y,\eta) \mapsto (y, r\eta)$ for any $r > 0$. 
The set $AZ$ is conic, 
and besides, it is invariant with respect to the mapping $(y,\eta) \mapsto (y, -\eta)$.

Denote by $n(A, \lambda)$ the number of singular values $s_n(A)$ 
greater than $\lambda>0$. 
The main result of the paper is the following Theorem.

\begin{T}\label{T}
	Let $\Omega\subset {\mathbb R}^d_+$, $\Sigma\subset\partial{\mathbb R}^d_+\times(0,\infty)$ be open bounded subsets satisfying
	$\overline\Omega\subset {\mathbb R}^d_+$, $\overline\Sigma\subset\partial{\mathbb R}^d_+\times(0,\infty)$, 
	$\mu_{d}(\partial\Sigma) = 0$ ($\mu_d$ is the $d$-dimensional Lebesgue measure on $\Sigma_0$, $\partial\Sigma$ is taken in the topology of $\Sigma_0$). Then
	\begin{equation}
		\lim_{\lambda \to 0} \lambda^{d} n({\cal O}_\Omega^\Sigma, \lambda) =
		(2\pi)^{-d} \iint_{AZ} dy d\eta \,\theta\left(\frac{1}{|\eta|\, |\eta_d|} - 1\right),
		\label{asymp}
	\end{equation}
	$\theta$ is the Heaviside function. 
\end{T}

Note that the integral in \re{asymp} is finite. Indeed, the Heaviside function vanishes if $\eta$ is large and $|\eta_d| \geqslant 1$; if $\eta$ is large and $|\eta_d| < 1$ then $(y,\eta) \notin AZ$ since $y\in\Omega$ and $\Omega$ is separated from $\partial{\mathbb R}^d_+$. We conclude that the Heaviside function is nonzero only on the bounded set of $(y,\eta)$.

The result~\re{asymp} implies the following behavior of singular values: 
\begin{equation}
	s_n({\cal O}_\Omega^\Sigma) \sim \left(\frac{\sigma}{n}\right)^{1/d}, \quad n\to\infty,
	\label{asymp_sn}
\end{equation}
where $\sigma$ is the r.h.s. of \re{asymp}.
Note that the asymptotics~\re{asymp_sn} has sense if $\sigma \ne 0$.
Also relation~\re{asymp} may be written in the equivalent form in terms of 
$({\cal O}_\Omega^\Sigma)^* {\cal O}_\Omega^\Sigma$:
\begin{equation}
	\lim_{\lambda \to 0} \lambda^{d/2} n(({\cal O}_\Omega^\Sigma)^* {\cal O}_\Omega^\Sigma, \lambda) = \sigma.
	\label{asymp_eig}
\end{equation}
Note that here $n$ 
coincides with the counting function of (positive) eigenvalues of $({\cal O}_\Omega^\Sigma)^* {\cal O}_\Omega^\Sigma$,
i.e. the number of eigenvalues greater than $\lambda$.
The relation~\re{asymp_eig} is a Weyl-type asymptotics for positive operator $({\cal O}_\Omega^\Sigma)^* {\cal O}_\Omega^\Sigma$.

In the rest of the paper we prove Theorem~\ref{T}.

\section{Operators $I$, $I_\pm$}
\label{sec_I}
Here we introduce and investigate FIOs connected with the problem~\re{BVP}.
Some basic concepts and facts of the theory of FIOs are used (canonical relations, composition of FIOs, principal symbol of composition of FIO with its adjoint)~-- we send the reader to books~\cite{Duist,Hor,Trev} for details.
First we express ${\cal O} v$ 
in terms of the following Cauchy problem for $\tilde u(x, t)$ in ${\mathbb R}^d\times {\mathbb R}$ (i.e. in the {\em whole} space):
\begin{align}
	&\partial^2_t \tilde u - \Delta \tilde u + \tilde q \tilde u = 0, \notag\\
	&\tilde u|_{t=0} = 0, \quad \partial_t \tilde u|_{t=0} = \tilde v.
	\label{Cauchy}
\end{align}
Here $\tilde q(x) = q(x', |x_d|)$ is an even continuation of $q$ to the whole space ($\tilde q$ is smooth as $q$ is compactly supported in ${\mathbb R}^d_+$), and $\tilde v \in C_0^\infty({\mathbb R}^d)$. 
Suppose $v\in C_0^\infty({\mathbb R}^d_+)$ and take $\tilde v$ such that $\tilde v = v$ in ${\mathbb R}^d_+$ and $\tilde v = 0$ in ${\mathbb R}^d\setminus{\mathbb R}^d_+$. Then the solutions of \re{BVP} and \re{Cauchy} are related as follows
\begin{equation}
	u(x, t) = \tilde u(x', x_d, t) + \tilde u(x', -x_d, t), \quad x_d > 0, \, t\in{\mathbb R}.
	\label{uu}
\end{equation}
Indeed,
$\tilde u(x', -x_d, t)$ satisfies the wave-type equation in \re{BVP} for $x_d > 0$ and has zero Cauchy data for $t=0$ in ${\mathbb R}^d_+$. The boundary condition for $x_d = 0$ is obviously satisfied.
It follows from \re{uu} that
\begin{equation}
	{\cal O} v = 2\, \tilde u|_{\Sigma_0}.
	\label{2u}
\end{equation}

Next we represent the solution of \re{Cauchy} in terms of FIOs: 
\[
	\tilde u = \tilde I_+ \tilde v + \tilde I_- \tilde v \textrm{ (modulo a smoothing operator)},
\]
\[
	(\tilde I_\pm \tilde v) (x, t) = (2\pi)^{-d} \iint_{{\mathbb R}^d\times{\mathbb R}^d} dy d\eta \, e^{i \tilde \varphi_\pm(x, t, y, \eta)}\, \tilde a_\pm(x, t, y, \eta)\, \tilde v(y),
\]
where $\tilde \varphi_\pm(x, t, y, \eta) = (x - y)\eta \pm t|\eta|$
is a nondegenerate phase function with phase variable $\eta$ and $\tilde a_\pm$ is an amplitude of order $-1$. 
The amplitude $\tilde a_\pm$ is classical:
\begin{equation}
	\tilde a_\pm \sim \sum_{j\geqslant 1} \zeta(|\eta|) \, \tilde a_\pm^{(j)},
	\label{ampl}
\end{equation}
$\zeta$ is a smooth function such that 
$\zeta(s) = 1$ for $s>2$ and $\zeta(s) = 0$ for $s < 1$,
functions $\tilde a_\pm^{(j)}$ are homogeneous of degree $-j$ in $\eta$.
Functions $\tilde a_\pm^{(j)}$ are solutions of certain transport equations~\cite{Maslov}, but we need only $\tilde a_\pm^{(1)}$:
\begin{equation}
	\tilde a_\pm^{(1)} = \pm (2i|\eta|)^{-1}. 
	\label{a1}
\end{equation}
Note that the definition of a FIO in \cite{Trev} implies the normalization factor $(2\pi)^{-d-1/4}$ rather than $(2\pi)^{-d}$ in formula for $\tilde I_\pm$.
Nevertheless our choice will be more convenient in further considerations.

The operators $\tilde I_\pm$ are continuous from $C_0^\infty({\mathbb R}^d)$ to $C^\infty({\mathbb R}^{d+1})$ (as well as from ${\cal E}'({\mathbb R}^d)$ to ${\cal D}'({\mathbb R}^{d+1})$).
Hence we may consider the restriction $\tilde I_\pm \tilde v|_{\Sigma_0}$ for $\tilde v\in C_0^\infty({\mathbb R}^d)$. 
The restriction can be written in the form of FIO as well:
\[
	(\check I_\pm \tilde v)(x', t) = (2\pi)^{-d} \iint_{{\mathbb R}^d\times{\mathbb R}^d} dy d\eta \, e^{i \varphi_\pm(x', t, y, \eta)}\, a_\pm(x', t, y, \eta)\, \tilde v(y), \quad (x',t) \in \Sigma_0,
\]
where $\varphi_\pm$ and $a_\pm$ are just the restrictions of $\tilde\varphi_\pm$ and $\tilde a_\pm$.
Note that
\[
	\varphi_\pm(x', t, y, \eta) = x' \eta' - y \eta \pm t|\eta|
\]
is a nondegenerate phase function.
The mapping
\[
	\tilde v \mapsto ( \tilde I_+ \tilde v + \tilde I_- \tilde v - \tilde u)|_{\Sigma_0} = \check I_+ \tilde v + \check I_- \tilde v - \tilde u|_{\Sigma_0}
\]
is smoothing. In view of \re{2u} we obtain the following representation of ${\cal O}$:
\begin{equation}
	{\cal O} = 2(\check I_+ + \check I_-)|_{C_0^\infty({\mathbb R}^d_+)} \textrm{ (modulo a smoothing operator)}.
	\label{OIcheck}
\end{equation}

Now let $\chi \in C_0^\infty(\Sigma_0)$ be a real-valued function, and put
\begin{equation}
	I_\pm = \chi\check I_\pm\big|_{C_0^\infty({\mathbb R}^d_+)}.
	\label{I}
\end{equation}
We consider $I_\pm$ as a FIO 
mapping functions in ${\mathbb R}^d_+$ to functions in $\Sigma_0$.
The phase function $\varphi_\pm$ and the amplitude $\chi a_\pm$ of the operator $I_\pm$ are functions of
$(x', t, y, \eta)$, where $(x', t)$ and $\eta$ vary over $\Sigma_0$ and ${\mathbb R}^d\setminus\{0\}$ correspondingly (just as for $\check I_\pm$), while $y$ varies over ${\mathbb R}^d_+$.

For the operator $I = 2(I_+ + I_-)$ due to~\re{OIcheck} and \re{I} we have
\begin{equation}
	\chi{\cal O} = I \textrm{ (modulo a smoothing operator)}.
	\label{OI}
\end{equation}
In sec.~\ref{secProof} some approximation of the characteristic function of the set $\Sigma$ will be chosen as the multiplier $\chi$ in~\re{I}. 
This will provide that the restriction of the operator $\chi{\cal O}$ to functions supported in $\Omega$ will approximate the operator ${\cal O}_\Omega^\Sigma$, which means that the analogous restriction of $I$ will approximate ${\cal O}_\Omega^\Sigma$.

The canonical relation of $I_\pm$ looks as follows
\[
	C_\pm = \{(x', t, \partial_{x',t}\varphi_\pm, y, -\partial_y\varphi_\pm) \,|\, \partial_\eta \varphi_\pm = 0 \} \subset
	(T^*\Sigma_0\setminus 0) \times (T^*{\mathbb R}^d_+\setminus 0).
\]
Note that $-\partial_y\varphi_\pm = \eta$, hence 
\[
	C_\pm = \{(x', t, \partial_{x',t}\varphi_\pm, y, \eta) \,|\, \partial_\eta \varphi_\pm = 0 \}
\]
and the phase variable $\eta$ may be treated as an element of $T^*_y{\mathbb R}^d_+$.
Introduce an open conic subset of $T^* {\mathbb R}^d_+\setminus 0$:
\[
	K = \{(y,\eta)\in T^* {\mathbb R}^d_+\setminus 0\,|\, \eta_d \ne 0\}.
\]
Consider the pair of canonical transformations from $K$ 
to $T^* \Sigma_0\setminus 0$:
\begin{equation}
	(y, \eta) \mapsto (x', t, \xi), \quad
	(x', t) = \gamma_\pm(y, \eta), \, \xi' = \eta', \, \xi_d = \pm|\eta|,
	\label{transform}
\end{equation}
where
\[
	\gamma_\pm(y, \eta) = \left(y' - y_d\, \frac{\eta'}{\eta_d},\, \pm y_d \, \frac{|\eta|}{\eta_d}\right).
\]
Analysis of $\varphi_\pm$ shows that the canonical relation $C_\pm$ is a graph of the mapping~\re{transform}.

Due to the multiplier $\chi$ in the definition~\re{I} we may consider the compositions $I_\pm^* I_\pm$ and $I_\mp^* I_\pm$. 
Since $C_\pm$ is a canonical graph the composition $I_\pm^* I_\pm$ is a $\Psi$DO (see~\cite{Trev}):
\[
	(I_\pm^* I_\pm v)(\overline y) = (2\pi)^{-d}\iint_{{\mathbb R}^d_+\times{\mathbb R}^d} dy d\eta\, e^{i (\overline y - y)\eta}\, b_\pm(\overline y, \eta) v(y).
\]
Since $I_\pm$ and $I_\pm^*$ are classical FIOs (i.e. Schwartz kernels can be represented by oscillatory integrals with classical amplitudes) the composition $I_\pm^* I_\pm$ is a classical $\Psi$DO:
\begin{equation}
	b_\pm \sim \sum_{j\geqslant 2} \zeta(|\eta|) \, b_\pm^{(j)},
	\label{bj}
\end{equation}
where $b_\pm^{(j)}$ are homogeneous of degree $-j$ in $\eta$.
The adjoint operator $I_\pm^*$ 
is a FIO with a canonical relation $C_\pm^t$.
Composition $I_\mp^* I_\pm$ is a smoothing operator in ${\mathbb R}^d_+$ since $C^t_\mp \circ C_\pm = \emptyset$.

Now we may write
\[
	I^* I = 4(I_+^* I_+ + I_-^* I_-) \textrm{ (modulo a smoothing operator)}
\]
and $I^* I$ is a classical $\Psi$DO in ${\mathbb R}^d_+$ of order $-2$ with symbol
\begin{equation}
	b = 4(b_++b_-).
	\label{bbb}
\end{equation}
From~\re{OI} it follows that
\begin{equation}
	(\chi{\cal O})^* \chi{\cal O} = I^* I 
	\textrm{ (modulo a smoothing operator)}.
	\label{hohohihi}
\end{equation}

We now calculate the principal symbol of the $\Psi$DO $I^*_\pm I_\pm$, i.e. the function $b_\pm^{(2)}(\overline y, \eta)$ in~\re{bj}.
We use the formula from \cite[sec. 8.6]{Trev} to express the principal symbol in terms of $\varphi_\pm$ and $a_\pm$. 
Consider the map
\[
	(x', t, y, \eta) \mapsto \left(y, -\partial_y\varphi_\pm(x', t, y, \eta), \partial_\eta\varphi_\pm(x', t, y, \eta)\right) =
	\left(y, \eta, \partial_\eta\varphi_\pm(x', t, y, \eta)\right)
\]
and denote by $\Delta_\pm$ the absolute value of its Jacobian determinant.
Then the principal symbol of $I^*_\pm I_\pm$ in $(\overline y,\eta)\in K$ equals
\[
	(|a_\pm^{(1)}|^2 \Delta_\pm^{-1}) (\gamma_\pm(\overline y, \eta), \overline y, \eta),
\]
where $a_\pm^{(1)}$ is the leading term of the amplitude $a_\pm$, due to~\re{a1} we have $a_\pm^{(1)} = \pm (2i|\eta|)^{-1}$.
Direct calculation shows that for $(\overline y, \eta)\in K$ we have $\Delta_\pm = |\eta_d|/|\eta|$, so
\[
	b_\pm^{(2)}(\overline y, \eta) = \frac{\chi(\gamma_\pm(\overline y, \eta))^2}{4 |\eta|\, |\eta_d|}.
\]
The principal symbol $b_\pm^{(2)}$ should be continued with zero to $(T^*{\mathbb R}^d_+\setminus 0) \setminus K$. 
Due to~\re{bbb} for the principal symbol $b^{(2)}$ in the expansion
\begin{equation}
	b \sim \sum_{j\geqslant 2} \zeta(|\eta|) \, b^{(j)}
	\label{bjj}
\end{equation}
we have
\begin{equation}
	b^{(2)}(\overline y, \eta) = \frac{\chi(\gamma_+(\overline y, \eta))^2 + \chi(\gamma_-(\overline y, \eta))^2}{|\eta|\, |\eta_d|}.
	\label{symbol_b}
\end{equation}

Further we will apply operators ${\cal O}$, $I$ to functions in $L_2(\Omega)$, which
is possible since these operators can be expressed in terms of FIOs that act on compactly supported distributions in ${\mathbb R}^d_+$ and $\Omega$ is separated from $\partial{\mathbb R}^d_+$.

\section{Spectral asymptotics of $\Psi$DOs}  
\label{sec_pdo}

We use the following result on spectral asymptotics of a $\Psi$DO
with symbol $\zeta(|\eta|) b^{(m)}(\overline y, \eta)$ 
\[
	(B v)(\overline y) = (2\pi)^{-d}\iint_{\Omega\times{\mathbb R}^d} dy d\eta\, e^{i (\overline y - y)\eta}\, \zeta(|\eta|)\, b^{(m)}(\overline y, \eta) v(y),
\]
where $b^{(m)}$ is homogeneous of degree $-m$ ($m>0$) in $\eta$ and $b^{(m)}\in C^\infty(\overline\Omega\times({\mathbb R}^d\setminus\{0\}))$, function $\zeta$ is the same as in~\re{ampl}. 
$B$ is a compact operator in $L_2(\Omega)$.
The result of~\cite{BS_anisotropic} claims that
\begin{equation}
	\lim_{\lambda\to 0} \lambda^{d/m} n(B, \lambda) = (2\pi)^{-d} \iint_{\Omega\times{\mathbb R}^d} d\overline y d\eta \,\theta(|b^{(m)}(\overline y, \eta)| - 1).
	\label{asymp_BS}
\end{equation}
Note that $B$ is not supposed to be self-adjoint or elliptic in~\cite{BS_anisotropic}.
Also note that \cite{BS_anisotropic} concerns the more general case 
when $b^{(m)}$ is anisotropically homogeneous and is only continuous in $\overline\Omega\times ({\mathbb R}^d\setminus\{0\})$.

Now let $I_\Omega$ be the restriction of $I$ to $L_2(\Omega)$.
The result~\re{asymp_BS} implies the following asymptotics for $I_\Omega^* I_\Omega$ 
(the operator $I_\Omega^* I_\Omega$ is compact since it is the restriction of the $\Psi$DO $I^* I$ of order $-2$ to $L_2(\Omega)$)
\begin{equation}
	\lim_{\lambda\to 0} \lambda^{d/2} n(I_\Omega^* I_\Omega, \lambda) = (2\pi)^{-d} \iint_{\Omega\times{\mathbb R}^d} d\overline y d\eta \,\theta(|b^{(2)}(\overline y, \eta)| - 1),
	\label{asymp_classic}
\end{equation}
where $b^{(2)}$ is given by~\re{symbol_b}.
To prove~\re{asymp_classic} we need to show that
\begin{equation}
	\lim_{\lambda\to 0} \lambda^{d/2} n(I_\Omega^* I_\Omega, \lambda) = \lim_{\lambda\to 0} \lambda^{d/2} n(B, \lambda),
	\label{IB}
\end{equation}
where $B$ is a $\Psi$DO in $\Omega$ with symbol $\zeta(|\eta|) b^{(2)}(\overline y, \eta)$.
The relation~\re{IB} claims merely 
that terms $b^{(j)}$, $j\geqslant 3$, in the expansion~\re{bjj}
do not influence the asymptotics of $n(I_\Omega^* I_\Omega, \lambda)$.
Although this fact seems to be trivial, the author could not find an appropriate reference, so a short proof is provided here.
We establish the estimate
\begin{equation}
	\lambda^{d/3} n(I_\Omega^* I_\Omega - B, \lambda) \leqslant C \quad \forall \lambda>0,
	\label{I-B}
\end{equation}
which means that $n(I_\Omega^* I_\Omega - B, \lambda)$ is estimated by less power of $1/\lambda$ than $(1/\lambda)^{d/2}$ and~\re{IB} will then follow
(see~\cite[sec. 11.6]{BS}).

Choose an open bounded set $\Omega'$ containing $\overline\Omega$, such that $\overline{\Omega'}\subset{\mathbb R}^d_+$, $\partial\Omega'\in C^\infty$, 
and denote by $I_{\Omega'}$ the restriction of $I$ to $L_2(\Omega')$.
Let $B'$ be the $\Psi$DO in $\Omega'$ with symbol $\zeta(|\eta|) b^{(2)}(\overline y, \eta)$.
Denote by $P_\Omega$ the projector onto $L_2(\Omega)$ acting in $L_2(\Omega')$.
We have
\begin{equation}
	n(I_\Omega^* I_\Omega - B, \lambda) = n(P_\Omega(I_{\Omega'}^*I_{\Omega'}-B')P_\Omega, \lambda). 
	\label{AAm}
\end{equation}
Choose $\rho \in C_0^\infty(\Omega')$, $\rho|_\Omega = 1$, and put $F = \rho (I_{\Omega'}^*I_{\Omega'}-B') \rho$. 
Since 
\[
	P_\Omega(I_{\Omega'}^*I_{\Omega'}-B')P_\Omega = P_\Omega F P_\Omega
\]
we have
\begin{equation}
	n(P_\Omega(I_{\Omega'}^*I_{\Omega'}-B')P_\Omega, \lambda) \leqslant n(F, \lambda).
	\label{AAm'}
\end{equation}
Here we used the following inequality
\[
	s_n(AL), s_n(LA) \leqslant \|L\| \, s_n(A)
\]
for a compact operator $A$ and a bounded operator $L$.
Put $D = (1\!\!\!\!1-\Delta)^{3/2}$, where $\Delta$ is the self-adjoint Laplace operator in $\Omega'$ with Dirichlet boundary condition. $D$ is a $\Psi$DO in $\Omega'$ \cite{Seeley}. 
The composition $F D$ is also well-defined as a $\Psi$DO of order $\leqslant 0$, hence it is a bounded operator in $L_2(\Omega')$. The operator $D^{-1}$ is also bounded in $L_2(\Omega')$.
For $v\in L_2(\Omega')$ we have
\[
	(Fv, Fv)_{L_2(\Omega')} = (FDD^{-1} v, FDD^{-1} v)_{L_2(\Omega')} \leqslant 
	\|F D\|^{2}_{L_2(\Omega')} (D^{-2} v, v)_{L_2(\Omega')}.
\]
Hence $F^*F \leqslant \|F D\|^{2}_{L_2(\Omega')} D^{-2}$ and
\[
	n(F^*F, \lambda) \leqslant n(\|F D\|^{2}_{L_2(\Omega')} D^{-2}, \lambda) = n(D^{-2}, \lambda/ \|F D\|^{2}_{L_2(\Omega')}).
\]
The well known asymptotics of eigenvalues 
of $-\Delta$ yields
\[
	n(D^{-2}, \lambda) \leqslant C \lambda^{-d/6}.
\]
We arrive at
\[
	n(F, \lambda) \leqslant C \lambda^{-d/3}.
\]
In view of \re{AAm}, \re{AAm'} this leads to~\re{I-B}.

\section{Proof of Theorem~\ref{T}}\label{secProof}
In this section we prove the relation~\re{asymp_eig}. 

First show that ${\cal O}^\Sigma_\Omega$ is a compact operator from $L_2(\Omega)$ to $L_2(\Sigma)$. 
Since the operator $I_\Omega^* I_\Omega$ is compact (see sec.~\ref{sec_pdo}) the operator $I_\Omega$ acting from $L_2(\Omega)$ to $L_2(\Sigma_0)$ is also compact.
Denote by ${\cal O}_\Omega$ the restriction of ${\cal O}$ to $L_2(\Omega)$.
Now due to~\re{OI} the composition $\chi{\cal O}_\Omega$ is a compact operator, 
which implies that ${\cal O}_\Omega^\Sigma$ is also compact since we can choose $\chi$ such that $\chi|_\Sigma = 1$.

To prove~\re{asymp_eig} we will make different choices of function $\chi$ in the definition~\re{I}.
First suppose that $0\leqslant\chi\leqslant 1$, $\chi = 1$ on $\Sigma$ and ${\rm supp}\,\chi\subset \partial{\mathbb R}^d_+\times (0,\infty)$.
For $v \in L_2(\Omega)$ we have
\[
	({\cal O}^\Sigma_\Omega v, {\cal O}^\Sigma_\Omega v)_{L_2(\Sigma)} \leqslant (\chi{\cal O}_\Omega\, v, \chi{\cal O}_\Omega\, v)_{L_2(\Sigma_0)}.
\]
Therefore $({\cal O}^\Sigma_\Omega)^* {\cal O}^\Sigma_\Omega \leqslant (\chi{\cal O}_\Omega)^* \chi{\cal O}_\Omega$ and so
\begin{equation}
	n(({\cal O}^\Sigma_\Omega)^* {\cal O}^\Sigma_\Omega, \lambda) \leqslant n((\chi{\cal O}_\Omega)^* \chi{\cal O}_\Omega, \lambda).
	\label{up}
\end{equation}
It follows that
\begin{equation}
	\varlimsup_{\lambda\to 0} \lambda^{d/2} n(({\cal O}^\Sigma_\Omega)^* {\cal O}^\Sigma_\Omega, \lambda) \leqslant \varlimsup_{\lambda\to 0} \lambda^{d/2} n((\chi{\cal O}_\Omega)^*\chi{\cal O}_\Omega, \lambda).
	\label{uplim}
\end{equation}
Due to~\re{hohohihi} the ``principal parts'' of operators $(\chi{\cal O}_\Omega)^*\chi{\cal O}_\Omega$ and $I_\Omega^* I_\Omega$ coincide, hence
\begin{equation}
	\lim_{\lambda\to 0} \lambda^{d/2} n((\chi{\cal O}_\Omega)^*\chi{\cal O}_\Omega, \lambda) = \lim_{\lambda\to 0} \lambda^{d/2} n(I_\Omega^* I_\Omega, \lambda) = \sigma_2,
	\label{sss}
\end{equation}
where $\sigma_2$ is the r.h.s. of~\re{asymp_classic}.
To prove this one should repeat the proof of~\re{IB}.

Now put
\[
	\kappa(y, \eta) = \chi(\gamma_+(y, \eta))^2 + \chi(\gamma_-(y, \eta))^2, \quad
	\kappa_\Sigma(y, \eta) = \chi_\Sigma(\gamma(y, \eta)),
\]
$\chi_\Sigma$ is a characteristic function of $\Sigma$.
We have
\begin{equation}
	\kappa_\Sigma \leqslant \kappa \leqslant 1.
	\label{ksik}
\end{equation}
Indeed, the first inequality is obvious if $\gamma(y, \eta)\notin\Sigma$.
Suppose $\gamma(y, \eta)\in\Sigma$, then $\gamma_+(y, \eta)\in\Sigma$ in case $\eta_d > 0$ and $\gamma_-(y, \eta)\in\Sigma$ in case $\eta_d < 0$. Now the inequality $\kappa_\Sigma \leqslant \kappa$ follows from $\chi|_\Sigma = 1$.
The second inequality follows from the fact that conditions $\gamma_+(y,\eta)\in{\rm supp}\,\chi$ and $\gamma_-(y,\eta)\in{\rm supp}\,\chi$ can not hold simultaneously since ${\rm supp}\,\chi\subset \partial{\mathbb R}^d_+\times (0,\infty)$.

Due to the first inequality in~\re{ksik} and formula~\re{symbol_b} we have $\sigma \leqslant \sigma_2$. Now we need to estimate $\sigma_2 - \sigma$.

Since $\overline\Omega\subset{\mathbb R}^d_+$ there exists $\varepsilon>0$ such that if $y\in\Omega$ and $|\eta_d|/|\eta| < \varepsilon$ then $\kappa(y, \eta) = 0$.
If $|\eta| > \varepsilon^{-1/2}$ and $y\in \Omega$ then 
$\kappa(y,\eta)/(|\eta| |\eta_d|) < 1$. Indeed, if $|\eta_d|/|\eta| \geqslant \varepsilon$ (otherwise $\kappa(y,\eta)=0$) then due to~\re{ksik}
\[
	\frac{\kappa(y,\eta)}{|\eta| \, |\eta_d|} \leqslant \frac{1}{|\eta| \, |\eta_d|} < 1.
\]
This means that the integrals in~\re{asymp_classic} and in~\re{asymp}
are taken over the set $y\in\Omega$, $|\eta| < \varepsilon^{-1/2}$.

The difference $\sigma_2-\sigma$ can be estimated by the measure of the set of $(y,\eta)$ such that
\begin{equation}
	|\eta| < \varepsilon^{-1/2}, \quad \kappa_\Sigma(y,\eta) < \kappa(y,\eta).
	\label{set}
\end{equation}
Due to~\re{ksik} the second condition implies that $\kappa_\Sigma(y,\eta) = 0$, and so $\gamma(y,\eta) \notin\Sigma$, which means that neither $\gamma_+(y,\eta)$ nor $\gamma_-(y,\eta)$ belongs to $\Sigma$. However, on the set~\re{set} we have $\kappa(y,\eta) > 0$, and so $\gamma_+(y,\eta)$ or $\gamma_-(y,\eta)$ belongs to ${\rm supp}\,\chi$.
We conclude that every $(y,\eta)$ from the set~\re{set} 
satisfies
\begin{equation}
	|\eta| < \varepsilon^{-1/2},\quad 
	\gamma_+(y,\eta) \textrm{ or } \gamma_-(y,\eta) \textrm{ belongs to } {\rm supp}\,\chi\setminus\Sigma.	
	\label{widerset}
\end{equation}
Now choosing $\chi$ such that ${\rm supp}\,\chi$ shrinks to $\overline\Sigma$ (note that in this case $\varepsilon$ does not depend on $\chi$) we make
the set~\re{widerset} shrink to the set of $(y,\eta)$ such that $|\eta| < \varepsilon^{-1/2}$ and $\gamma_+(y,\eta)$ or $\gamma_-(y,\eta)$ belongs to $\partial\Sigma$.
This set has zero measure in $\Omega\times{\mathbb R}^d$ since $\mu_d(\partial\Sigma) = 0$.
This means that choosing appropriate $\chi$ we can provide that 
the difference $\sigma_2-\sigma$ is arbitrarily small.
Together with~\re{uplim} and \re{sss} this means that
\begin{equation}
	\varlimsup_{\lambda\to 0} \lambda^{d/2} n(({\cal O}^\Sigma_\Omega)^* {\cal O}^\Sigma_\Omega, \lambda) \leqslant \sigma.
		\label{upp}
\end{equation}

Taking $\chi\in C_0^\infty(\Sigma)$, $0\leqslant\chi\leqslant 1$, we obtain the inequality reverse to~\re{up}, therefore we have the relation (which is a counterpart of~\re{uplim})
\[
	\varliminf_{\lambda\to 0} \lambda^{d/2} n(({\cal O}^\Sigma_\Omega)^* {\cal O}^\Sigma_\Omega, \lambda) \geqslant \varliminf_{\lambda\to 0} \lambda^{d/2} n((\chi{\cal O}_\Omega)^*\chi{\cal O}_\Omega, \lambda).
\]
Then arguing the same way as in proof of~\re{upp} we obtain that
\[
	\varliminf_{\lambda\to 0} \lambda^{d/2} n(({\cal O}^\Sigma_\Omega)^* {\cal O}^\Sigma_\Omega, \lambda) \geqslant \sigma.
\]
Together with~\re{upp} this yields~\re{asymp_eig} and thus Theorem~\ref{T} is proved.

Now we expose the reason why we used the smooth multiplier $\chi$ in the definition~\re{I} of the operator $I_\Omega$. 
Instead of this we could try to deal with $I'_\Omega := 2(\check I_+ + \check I_-)$ considered as an operator from $L_2(\Omega)$ to $L_2(\Sigma)$ (see formula~\re{OIcheck}).
It seems to be possible 
to 
consider $(I'_\Omega)^* I'_\Omega$
as a $\Psi$DO with the principal symbol
\[
	\frac{\chi_\Sigma(\gamma(\overline y, \eta))}{|\eta|\, |\eta_d|}
\]
($\chi_\Sigma$ is the characteristic function of $\Sigma$), which is {\em discontinuous} in $\Omega\times({\mathbb R}^d\setminus\{0\})$.
However, in this case the result~\re{asymp_BS} of~\cite{BS_anisotropic} can not be applied since the principal symbol is required to be continuous there.
Note that classical results on spectral asymptotics of $\Psi$DOs with discontinuous symbols by H.~Widom and their improvements also can not be applied to $(I'_\Omega)^* I'_\Omega$~-- we do not go into details here.

\end{document}